\documentclass{amsart}
\usepackage{hyperref}
\usepackage{graphicx,amsmath,amssymb,amsthm}
\usepackage{color}
\def\argmax{\operatornamewithlimits{arg\,max}}
\def\<#1,#2>{\langle#1,#2\rangle}
\newcommand{\appr}{^{^{_\sim}}}
\newcommand{\apprd}{^{^{_{\sim\sim}}}}
\newcommand{\mrm}[1]{\text{\rm #1}}
\newcommand{\dd}{\,d}
\newcommand{\col}{\operatorname{col}}

\newcommand{\new}[1]{{\em #1}}

\newcommand{\R}{\mathbb{R}}
\newcommand{\rbar}{\overline{\R}}
\newcommand{\Z}{\mathbb{Z}}
\newcommand{\sK}{\mathcal{K}}

\newcommand{\im}{\mathrm{im}\,}

\newcommand{\projimker}[2]{\Pi_{#1}^{#2}}
\newcommand{\rmax}{\mathbb{R}_{\max}}
\newcommand{\rmaxb}{\overline{\R}_{\max}}
\newcommand{\rminb}{\overline{\R}_{\min}}

\newcommand{\supp}{\mathop{\text{\Large$\vee$}}}
\newcommand{\inff}{\mathop{\text{\Large$\wedge$}}}
\newcommand{\maxx}{\max}

\newcommand{\comp}{\circ}

\newcommand{\bydef}{\stackrel{\mathrm{def}}{=}}
\newcommand{\set}[2]{\{#1\mid\,#2\}}
\newcommand{\lres}{\backslash}

\newcommand{\sh}{^{\sharp}}

\newcommand{\calu}{\mathcal{U}}
\newcommand{\calv}{\mathcal{V}}
\newcommand{\calw}{\mathcal{W}}
\newcommand{\calx}{\mathcal{X}}
\newcommand{\caly}{\mathcal{Y}}
\newcommand{\calz}{\mathcal{Z}}
\def\<#1,#2>{\langle#1\mid #2\rangle}
\newtheorem{thm}{Theorem}
\newtheorem{prop}[thm]{Proposition}
\newtheorem{lem}[thm]{Lemma}

\theoremstyle{definition}
\newtheorem{exmp}[thm]{Example}

\newtheorem{defin}[thm]{Definition}
\title[A max-plus finite element method]{A max-plus finite element method for solving finite horizon deterministic optimal control problems}
\author{Marianne Akian}
\author{St\'ephane Gaubert}
\author{Asma Lakhoua}
\address{INRIA, Domaine de Voluceau, 78153 Le Chesnay C\'edex, France}
\email{\{Marianne.Akian,Stephane.Gaubert,Asma.Lakhoua\}@inria.fr}
\date{March 15, 2004. Revised April 7, 2004. Prepared for MTNS'04.}

\begin{document}
\begin{abstract}
We introduce a max-plus analogue of the Petrov-Galerkin finite
element method, to solve finite horizon deterministic
optimal control problems. The method relies on a 
max-plus variational formulation, and exploits the properties of projectors
on max-plus semimodules. 
We obtain a nonlinear discretized semigroup,
corresponding to a zero-sum two players game.
We give an error estimate of order $\sqrt{\Delta t}+\Delta x(\Delta t)^{-1}$,
for a subclass of problems in dimension 1. 
We compare our method with a max-plus based discretization method
previously introduced by Fleming and McEneaney.
\end{abstract}
\maketitle
\section{Introduction}
We consider the optimal control problem:
\begin{subequations}
\label{problemP}
\begin{align}
\label{p1}
\mrm{maximize }
\int_0^T \ell(x(s),u(s))\dd s+\phi(x(T))
\end{align}
over the set of trajectories $(x(\cdot),u(\cdot))$ satisfying
\begin{align}
\label{p2}
\dot{x}(s)&=f(x(s),u(s)),\quad
x(0)=x,\quad x(s)\in X,\quad u(s)\in U \enspace,
\end{align}
\end{subequations}
for all $0\leq s\leq T$.
Here, the \new{state space} $X$ is a subset of $\R^n$,
the set of \new{control values} $U$ is a subset of $\R^m$,
the \new{horizon} $T>0$ and the \new{initial condition} $x\in X$ are given,
we assume that the map $u(\cdot)$ is measurable, 
and that the map $x(\cdot)$ is absolutely continuous.
We also assume that the \new{instantaneous reward}
or \new{Lagrangian} $\ell:X\times U \to \R$,
and the \new{dynamics} $f:X\times U \to \R^n$, are
sufficiently regular maps, and that the \new{terminal reward}
$\phi$ is a map $X\to \R\cup\{-\infty\}$.
The \new{value function} $v$ associates to any $(x,t)\in X\times [0,T]$
the supremum $v(x,t)$ of $\int_0^t \ell(x(s),u(s))\dd s+\phi(x(t))$,
under the constraint~\eqref{p2}, for $0\leq s \leq t$.
Under certain regularity assumptions, it is known that $v$ is solution of the
Hamilton-Jacobi equation
\begin{subequations}\label{HJ}
\begin{gather}
-\frac{\partial v}{\partial t}+H(x,\frac{\partial v}{\partial x})=0, \quad (x,t)
\in X \times (0,T] \enspace,\label{HJ1}\end{gather}
with initial condition:
\begin{gather}
v(x,0)=\phi(x), \quad  x \in X \enspace,
\end{gather}
\end{subequations}
where $H(x,p)=\sup_{u\in U}\ell(x,u)+p\cdot f(x,u)$ is the \new{Hamiltonian}
of the problem (see for instance~\cite{lions,soner,barles}).
The \new{evolution semigroup} $S^t$ of~\eqref{HJ} associates
to any map $\phi$ the function $v^t:=v(\cdot,t)$, where $v$ is the value
function of the optimal control problem~\eqref{p1}. 

Maslov~\cite{maslov73} (see also~\cite{maslov92,kolokoltsov})
observed that the evolution semigroup $S^t$ is max-plus linear.
Recall that the \new{max-plus semiring}, $\rmax$, is the set
$\R\cup\{-\infty\}$, equipped
with the addition $a\oplus b=\max(a,b)$ and the multiplication
$a\otimes b=a+b$. By \new{max-plus linearity},
we mean that for all maps $f,g$ from $X$ to $\rmax$, 
and for all $\lambda\in\rmax$, we have 
\begin{align*}
S^t(f\oplus g)&=S^tf\oplus S^tg \enspace ,\\
S^t(\lambda f)&=\lambda S^tf \enspace ,
\end{align*}
where $f\oplus g$ denotes the map $x\mapsto f(x)\oplus g(x)$,
and $\lambda f$ denotes the map $x\mapsto \lambda \otimes f(x)$.
Linear operators over max-plus type semirings have been widely studied,
see for instance~\cite{cuning,maslov92,baccelli,kolokoltsov,gondran-minoux}.

In this paper, we introduce a new discretization method to solve the
deterministic optimal control problem~\eqref{problemP},
using the max-plus linearity of the semigroup $S^t$.
In~\cite{mceneaney}, Fleming and McEneaney introduced
a max-plus based discretization method to solve a subclass
of Hamilton-Jacobi equations
(with a Lagrangian $\ell$ quadratic with respect to $u$, and a dynamics
$f$ affine with respect to $u$). They approximated the evolution
semigroup $S^t$ by a max-plus linear semigroup acting on a
finitely generated semimodule of functions. This work
was pursued in~\cite{mceneaney00,mceneaney00b,mceneaney02,mceneaney03}. 
Another max-plus based numerical work on Hamilton-Jacobi equations
is due to Bacaer~\cite{bacaer01,bacaer02}.
The different discretization that we introduce here
relies on a notion of max-plus ``variational formulation'',
which originates from the notion
of generalized solution of Hamilton-Jacobi
equations of Maslov and Kolokoltsov~\cite{kolokltsovmaslov88},
\cite[Section 3.2]{kolokoltsov}. This discretization,
which can be interpreted geometrically in terms
of projections on semimodules,
is similar to the classical finite element method.
We shall see that the space of test functions
must be different from the space in which
the solution is represented, so that our discretization
is indeed a max-plus analogue of the
Petrov-Galerkin finite element method. We illustrate
the method by numerical examples. We also give an error
estimate, in dimension one, of order
$\sqrt{\Delta t}+\Delta x(\Delta t)^{-1}$,
which is the same as the order obtained
for existing discretization methods, see~\cite{falcone} and~\cite[Appendix A, by M. Falcone]{bardi-capuzzo-dolcetta} 

The present paper is only a preliminary account:
the results will be detailed elsewhere. A first
presentation of the method appeared in~\cite{asma}.
\section{Preliminaries on residuation and projections over semimodules}
In this section we recall some classical residuation results (see
for example \cite{marie-louise}, \cite{birkhoff}, \cite{blyth},
\cite{baccelli}), and their application to linear maps on idempotents
semimodules (see~\cite{litvinov,ilade}). We also review some
results of \cite{wodes,ilade} concerning projectors over
semimodules.
\subsection{Residuation, semimodules, and linear maps}
If $(S,\leq)$ and $(T,\leq)$ are (partially) ordered sets, we say that a map
$f:S\to T$ is \new{monotone} if $s \leq s' \implies f(s)\leq
f(s')$. We say that $f$ is \new{residuated} if there exists a map
$f\sh: T\to S$ such that
\[
f(s) \leq t \iff s\leq f\sh(t) \enspace .
\]
The map $f$ is residuated if, and only if, for all $t\in T$,
$\set{s\in S}{f(s)\leq t}$ has a maximum element in $S$. Then,
 \begin{align*}
f\sh(t)&=\max\set{s\in S}{f(s)\leq t},\quad \forall t\in T \enspace.
\end{align*}
Moreover, in that case, we have 
\begin{align}\label{ffdf}
&f\comp f\sh\comp f=f\sh\text{ and } f\sh\comp f\comp f\sh=f\enspace.
\end{align}
If a set $\sK$ is a monoid for a commutative idempotent law
$\oplus$ (\new{idempotent} means that $a\oplus a =a$),
the \new{natural order} on $\sK$ is defined
by $a\leq b \iff a\oplus b=b$. We say that $\sK$ is \new{complete}
as a naturally ordered set 
if any subset of $\sK$ has a least upper bound for the natural
order. If $(\sK,\oplus,\otimes)$ is an idempotent semiring,
i.e., a semiring whose addition is idempotent, we say that
the semiring $\sK$ is \new{complete} if it is complete
as a naturally ordered set, and if the left and right multiplications,
$L^\sK_a$, $R^\sK_a: \sK\to \sK$, $L^\sK_a(x)=ax$, $R^\sK_a(x)=xa$,
are residuated.

The max-plus semiring, $\rmax$,
is an idempotent semiring. It is not complete, but it can be embedded
in the complete idempotent semiring  $\rmaxb$ obtained 
by adjoining $+\infty$ to $\rmax$, with the convention that
$-\infty$ is absorbing for the multiplication $a\otimes b=a+b$.
The map $x\mapsto -x$ from $\rbar$ to itself yields an isomorphism
from $\rmaxb$ to the complete idempotent semiring $\rminb$,
obtained by replacing $\max$ by $\min$ and by exchanging the roles
of $+\infty$ and $-\infty$ in the definition of $\rmaxb$.

Semimodules over semirings are defined like modules over rings,
mutatis mutandis, see~\cite{litvinov,ilade}.
When $\sK$ is a complete idempotent semiring, we say that a (right)
$\sK$-semimodule $\calx$ is \new{complete} if it is complete as a
naturally ordered set, and if, for all $u\in \calx$ and $\lambda\in \sK$,
the right and left multiplications, $R^\calx_{\lambda}:\;\calx\to \calx$,
$v\mapsto v\lambda$ and 
$L^\calx_{u}:\;\sK\to \calx$, $\mu\mapsto u\mu$, are residuated.
In a complete semimodule $\calx$, we define, for all $u,v\in \calx$,
\begin{align*}
  u\lres v &\bydef (L_u^\calx)\sh(v) = \max\set{\lambda\in \sK}{u\lambda \leq v} \enspace .
\end{align*}

We shall use \new{semimodules of functions}:
when $X$ is a set and $(\sK,\oplus,\otimes)$
is a complete idempotent semiring, the set of functions
$\sK^X$ is a complete $\sK$-semimodule for the 
componentwise addition $(u,v)\mapsto u\oplus v$ (defined by 
$(u\oplus v)(x)= u(x)\oplus v(x)$),
and the componentwise multiplication 
$(\lambda,u)\mapsto u \lambda$ (defined by $(u\lambda)(x)= u(x)\otimes
\lambda$).

If $\sK$ is an idempotent semiring, 
and if $\calx$ and $\caly$ are $\sK$-semimodules,
we say that a map $A:\calx\to \caly$ is \new{additive} if for all $u,v\in \calx$,
$A(u\oplus v)=A(u)\oplus A(v)$ and that $A$
 is \new{homogeneous} if for all $u\in \calx$ and $\lambda\in \sK$,
$A (u \lambda)= A(u)\lambda$.
We say that $A$ is \new{linear},
or is a \new{linear operator}, if it is additive and homogeneous.
Then, as in classical algebra, we use the
notation $Au$ instead of $A(u)$. 
When $A$ is residuated and
$v\in \caly$, we use the notation $A\backslash v$ or $A\sh v$ instead of
$A\sh (v)$.
We denote by $L(\calx,\caly)$ the set of linear operators
from $\calx$ to $\caly$.
If $\sK$ is a complete idempotent semiring,
if $\calx,\caly,\calz$ are complete $\sK$-semimodules,
and if $A\in L(\caly,\calz)$ is residuated, then 
the map $L_A:L(\calx,\caly)\to L(\calx,\calz),\; B\mapsto A\comp B$, 
is residuated and we set $A\backslash C:=(L_A)\sh(C)$,
for all $C\in L(\calx,\calz)$.

If $X$ and $Y$ are two sets, $\sK$
is a complete idempotent semiring, and $a\in \sK^{X\times Y}$,
we construct the linear operator $A$ from $\sK^Y$ to $\sK^X$ which associates
to any $u\in \sK^Y$ the function $Au\in \sK^X$ such that
$Au(x)=\supp_{y\in Y} a(x,y)\otimes u(y)$, where $\vee$ denotes
the supremum for the natural order.
We say that $A$ is the \new{kernel operator} with \new{kernel}
or \new{matrix} $a$.
We shall often use the same notation $A$ for the operator and the kernel.
As is well known (see for instance~\cite{baccelli}), the kernel
operator $A$ is residuated, and
\[
(A\backslash v)(y)=\inff_{x\in X}A(x,y)\backslash v(x),
\]
where $\wedge$ denotes the infimum for the natural order.
In particular, when $\sK=\rmaxb$, we have
\begin{align}
\label{e-conv}
(A\backslash v)(y)=\inff_{x\in X}(-A(x,y)+ v(x))= [- A^* (-v)](y)
\enspace,
\end{align}
where $A^*$ denotes the \new{transposed operator} $\sK^X\to \sK^Y$, 
which is associated to the kernel $A^*(y,x)=A(x,y)$. 
(In~\eqref{e-conv}, we use the convention that $+\infty$ is
absorbing for addition.)

\subsection{Projectors on semimodules}
Let $\calv$ denote a \new{complete subsemimodule} of a complete semimodule
$\calx$ over a complete idempotent semiring $\sK$, i.e., a subset of $\calx$
that is stable by arbitrary sups and by the action of scalars.
We call \new{canonical projector} on $\calv$ the map
\begin{equation}\label{projecteur}
P_\calv: \calx\to \calx, \quad u\mapsto 
P_\calv(u) = \maxx\set{v\in \calv}{v\leq u}. 
\end{equation}
Let $W$ denote a \new{generating family} of a complete
subsemimodule $\calv$, which means that
any element $v\in \calv$ can be written as
$v=\supp\set{w\lambda_w}{w\in W}$, for some $\lambda_w\in\sK$.
It is known that
\[
P_\calv(u) =  \supp_{w\in W} w (w\lres u) 
\]
(see for instance~\cite{ilade}).
If $B:\calu\to\calx$ is a residuated linear operator, then the image $\im B$
of $B$ is a complete subsemimodule of $\calx$, and
\begin{equation}\label{PimB}
P_{\im B}=B\comp B\sh.
\end{equation}
The max-plus finite element methods relies
on the notion of projection on an image, parallel
to a kernel, which was introduced by Cohen, 
the second author, and Quadrat, in~\cite{wodes}.
The following theorem, of which Proposition~\ref{vhdelta} below
is an immediate corollary, is a variation
on the results of~\cite[Section~6]{wodes}.
\begin{thm}[Projection on an image parallel to a kernel]\label{piBC}
Let $B:\calu\to\calx$ and $C:\calx\to\caly$ be two residuated linear
operators. Let $\projimker{B}{C}=B\comp(C\comp B)\sh\comp C$. We have
$\projimker{B}{C}=\projimker{B}{}\comp \projimker{}{C}$, where
$\projimker{B}{}=B\comp B\sh$ and $\projimker{}{C}=C\sh\comp
C$. Moreover, $\projimker{B}{C}$ is a projector
$\big((\projimker{B}{C})^2=\projimker{B}{C}\big)$, and for all
$x\in\calx$:
\[
\projimker{B}{C}(x)=\maxx\set{y\in\im B}{Cy\leq Cx}.
\]
\end{thm}
The results of~\cite{wodes} characterize
the existence and uniqueness, for all $x\in X$,
of $y\in \im B$ such that $Cy=Cx$. In that case,
$y=\projimker{B}{C}(x)$.

When $\sK=\rmaxb$, and $C:\rmaxb^X\to\rmaxb^Y$
is a kernel operator,
$\projimker{}C=C\sh\comp C$ has an interpretation similar
to~\eqref{PimB}:
\[ 
\projimker{}{C}(v)=C\sh\comp C(v)=-P_{\im C^*}(-v)
=P_{-\im C^*} (v)\enspace,\]
where $-\im C^*$ is thought of as a $\rminb$-subsemimodule
of $\rminb^X$, so that, 
\[
P_{-\im C^*} (v)= \min\set{w\in -\im C^*}{w\geq v} \enspace .
\]
where $\leq$ denotes here the usual order
on $\rbar^X$, since the natural order of 
$\rminb^X$ is the reverse of the usual order.
When $B:\rmaxb^U\to\rmaxb^X$ is also a kernel
operator, we have
\[
\projimker BC=P_{\im B}\comp P_{-\im C^*} \enspace .
\]
This factorization will be instrumental in the geometrical
interpretation of the finite element algorithm, see Example~\ref{ex-dist}
below.
\section{The max-plus finite element method}
\subsection{Max-plus variational formulation}
We now describe the max-plus finite element method to solve
the optimal control problem~\eqref{p1}.
Let $S^t$ and $v^t$ be defined as in the introduction.
Using the semigroup property $S^{t+t'}=S^t\circ S^{t'}$, for $t,t'>0$, 
we have the recursive equation:
\begin{equation}\label{exact}
\begin{array}{cc}
v^{t+\Delta t}=S^{\Delta t} v^t, & t=0,\Delta t, \cdots, T-\Delta t
\end{array}
\end{equation}
with $v^0=\phi$ and $\Delta t=\frac{T}{N}$, for some
positive integer $N$. Let $\calw$ be
a $\rmaxb$-semimodule of functions from $X$ to $\rmaxb$ such that
$\phi\in\calw$ and for all $v\in\calw$, $t>0$, $S^tv\in\calw$. We
suppose given a ``dual'' semimodule $\calz$ of ``test functions'' from $X$
to $\rmaxb$. The max-plus \new{scalar product}
is defined by $\<u,v>=\sup_{x\in X} u(x)+v(x)$,
for all functions $u,v:X\to \rbar$, with the convention that $-\infty$
is absorbing for the addition $+$.
We replace \eqref{exact} by:
\begin{equation}\label{produitscalaire}
\<z,v^{t+\Delta t}>=\<z,S^{\Delta t}v^{t}> ,
\quad  \forall z\in\calz,\quad 
t=0,\Delta t,\ldots,T-\Delta t \enspace,
\end{equation}
with $v^{\Delta t},\ldots,v^T\in\calw$.
Equation~\eqref{produitscalaire}
can be seen as the analogue of a \new{variational}
or\new{weak formulation}. Kolokoltsov and Maslov used this formulation
in \cite{kolokltsovmaslov88} and \cite[section 3.2]{kolokoltsov}
to define a notion of generalized solution of Hamilton-Jacobi equations.
\subsection{Ideal max-plus finite element method}
We consider a semimodule $\calw_h\subset\calw$ with
generating family $\{w_i\}_{1\leq i\leq p}$. We call \new{finite elements}
the functions $w_i$. We approximate $v^t$ by
$v_h^t\in\calw_h$, that is:
\[
v^t\simeq v_h^t=\supp_{1\leq i \leq p}w_i \lambda^t_i\enspace,
\]
where $\lambda_i^t\in\rmax$. We also consider a semimodule
$\calz_h\subset\calz$ with generating family $\{z_j\}_{1\leq j\leq
  q}$. The functions $z_1,\cdots,z_q$ will act as test functions. We replace \eqref{produitscalaire} by
\begin{equation}\label{pdtscalapproch}
\begin{array}{cc}
\<z_j,v_h^{t+\Delta t}>=\<z_j,S^{\Delta t}v_h^{t}>, & \forall 1\leq
j\leq q\enspace,\\
\end{array}
\end{equation}
for $t=0,\Delta t, \cdots, T-\Delta t$, with $v_h^0=\phi_h\simeq\phi$ and
$v_h^t\in\calw_h$, $t=0,\Delta t,\cdots,T$.

Since Equation \eqref{pdtscalapproch} need not have a
solution, we look for the maximal subsolution,
i.e.\ the maximal solution $v_h^{t+\Delta t}\in \calw_h$ of
\begin{subequations}\label{inegal}
\begin{align}
\<z_j,v_h^{t+\Delta t}> \quad \leq \quad \<z_j,S^{\Delta t}v_h^{t}> \quad 
\forall 1\leq j\leq q\enspace. \label{infouegal}
\end{align}
We also take for the approximate value function $v_h^0$ at time $0$
the maximal solution $v_h^0\in \calw_h$ of
\begin{align}
v_h^{0}\leq v^0\enspace.\label{inegvh0}
\end{align}
\end{subequations}
Let us denote by $W_h$ the max-plus linear operator from
$\rmax^p$ to $\calw$ with matrix $W_h=\col(w_{i})_{1\leq i\leq p}$,
and by $Z_h^*$ the max-plus linear operator from $\calw$ to $\rmaxb^q$ 
whose transposed matrix is $Z_h=\col(z_j)_{1\leq j\leq q}$.
This means that 
$W_h\lambda=\supp_{1\leq i \leq p}w_i \lambda_i$ for all
$\lambda=(\lambda_i)_{i=1,\ldots, p}\in \rmax^p$, and
$(Z_h^* v)_j=\< z_j,v>$ for all $v\in\calw$ and $j=1,\ldots,q$.
Applying Theorem~\ref{piBC} to $B=W_h$ and $C=Z_h^*$ and using
$\calw_h=\im W_h$, we get:
\begin{prop}\label{vhdelta}
The maximal solution $v_h^{t+\Delta t}\in \calw_h$
of~\eqref{infouegal} is given by 
$v_h^{t+\Delta t}=S_{h}^{\Delta t}v_h^t$, where 
\[
S_{h}^{\Delta t}=\projimker{W_h}{Z_h^*}\comp 
 S^{\Delta t}\enspace.
\]
\end{prop}
\begin{prop}Let $v_h^{t}\in \calw_h$ be the  maximal solution 
of~\eqref{inegal}, for $t=0,\Delta t,\ldots, T$. 
Then,  for every $t=0,\Delta t,\ldots, T$,
there exists $\lambda^t\in\rmax^p$ such that
$v_h^t=W_h\lambda^t$. Moreover, the maximal $\lambda^t$ satisfying these
conditions verifies the recursive equation
\begin{subequations}
\begin{align}\label{lambdat+dt}
\lambda^{t+\Delta t}=(Z_h^*W_h)\backslash (Z_h^*S^{\Delta t}W_h\lambda^t)
\enspace,
\end{align}
with the initial condition:
\begin{align*}
\lambda^{0}=W_h\backslash \phi\enspace.
\end{align*}
\end{subequations}
\end{prop}
\begin{proof}
Since $v_h^t\in \calw_h$, $v_h^t=W_h\lambda^t$,
and the maximal $\lambda^t$ satisfying this condition is
$\lambda^t=W_h\sh (v_h^t)$, for all $t=0,\Delta t,\ldots, T$.
Since  $v_h^0$ is the maximal solution of~\eqref{inegvh0},
then by~\eqref{projecteur} 
and~\eqref{PimB}, $v_h^0=P_{\calw_h}(\phi)=W_h \comp W_h\sh (\phi)$, hence
$\lambda^0= W_h\sh\comp W_h\comp W_h\sh (\phi)=W_h\sh(\phi)$. 
Let $t=0,\ldots, T_\Delta t$. Using Proposition~\ref{vhdelta},  
Theorem~\ref{piBC},~\eqref{ffdf} and the property that
$(f\comp g)\sh=g\sh\comp f\sh$ for all residuated maps $f$ and $g$,
we get 
\begin{eqnarray*}
\lambda^{t+\Delta t} & = & W_h\sh \comp \projimker{W_h}{Z_h^*}
\comp S^{\Delta t}(W_h\lambda^t )\\
& = & W_h\sh\comp W_h\comp  W_h\sh \comp (Z_h^*)\sh\comp Z_h^*
\comp S^{\Delta t}(W_h\lambda^t)\\
&= & W_h\sh \comp (Z_h^*)\sh\comp Z_h^*\comp S^{\Delta t} (W_h\lambda^t)\\
& = & (Z_h^* W_h)\sh (Z_h^* S^{\Delta t}W_h\lambda^t)\enspace.
\end{eqnarray*}
which yields~\eqref{lambdat+dt}.
\end{proof}

The maps $A_h:=Z_h^*W_h:\rmax^p\to \rmax^q$ 
and $B_{h}:=Z_h^* S^{\Delta t}W_h:\rmax^p\to \rmax^q$
are max-plus linear operators,
and the entries of their corresponding matrices 
are given, for $1\leq i \leq p$ and $1\leq j\leq q$, by:
\begin{align}
(A_h)_{ji}&=\<z_j,w_i> \label{matrixA}\\
(B_{h})_{ji}&=\<z_j,S^{\Delta   t}w_i> \label{matrixB}
\\
&=\<(S^*)^{\Delta t} z_j, w_i>
\enspace,\label{matrixB'}
\end{align}
where $S^*$ is the \new{transposed semigroup} of $S$, which
is the evolution semigroup associated to the optimal
control problem in which the sign of the dynamics is changed.

The ideal max-plus finite element method can be summarized as follows:
\begin{enumerate}
\item Choose $\Delta t=\frac{T}{N}$ and
the finite elements $(w_i)_{1\leq i\leq p}$ and $(z_j)_{1\leq j\leq q}$, 
\item Compute the matrix $A_h$ by~\eqref{matrixA}
and the matrix $B_{h}$ by~\eqref{matrixB} or by~\eqref{matrixB'},
\item Compute $\lambda^{0}=W_h\backslash\phi$ and
 $v_h^{0}=W_h\lambda^{0}$.
\item For $t=\Delta t, 2\Delta t,\ldots,T$,
compute $\lambda^{t}=A_h\backslash (B_{h}\lambda^{t-\Delta t})$ and
 $v_h^{t}=W_h\lambda^{t}$.
\end{enumerate}
Then, $v_h^t$ approximates the value function at time $t$, $v^t$.

The recursion $\lambda^{t}=A_h\backslash (B_{h}\lambda^{t-\Delta t})$
may be written explicitly as 
\[
\lambda^t_i = \min_{1\leq j\leq q}
\Big( -(A_{h})_{ji}+ \max_{1\leq k\leq p}\big( (B_h)_{jk} + \lambda^{t-\Delta t}_k\big) \Big),
\quad \mrm{for } 1\leq i\leq p\enspace .
\]
Observe that this recursion 
may be interpreted as the dynamic programming equation of a deterministic
zero-sum two players game, with finite action and state spaces.

In order to implement this method, we must specify how to
compute the entries of $A_h$ and $B_h$ in~\eqref{matrixA}
and~\eqref{matrixB} or~\eqref{matrixB'}.
In some cases, these computations can be done analytically.
Computing $A_h$ from~\eqref{matrixA} is an optimization
problem which may be solved by  standard algorithms.
We shall discuss in the following section the approximation
of $B_h$. 
\subsection{Effective max-plus finite element method}
We first discuss the approximation of $S^{\Delta t}w$ 
for every finite element $w$.
The Hamilton-Jacobi equation~\eqref{HJ1} suggests to
approximate $S^{\Delta t}w$ by the function $[S^{\Delta t}w]\appr $ such that
\begin{equation}\label{stilde}
[S^{\Delta t}w]\appr (x)=w(x)+\Delta t H(x,\frac{\partial{w}}{\partial x}),
\quad \mrm{for all} x\in X .
\end{equation}
Let $[S^{\Delta t}W_h]\appr $ denotes the max-plus linear operator from
$\rmax^p$ to $\calw$ with matrix $[S^{\Delta t}W_h]\appr =
\col([S^{\Delta t}w_i]\appr )_{1\leq i\leq p}$, which means that
\[ 
[S^{\Delta t}W_h]\appr \lambda=\supp_{1\leq i \leq p}[S^{\Delta t}w_i]\appr
\lambda_i
\]
for all $\lambda=(\lambda_i)_{1\leq i\leq p}\in \rmax^p$.
The above approximation of  $S^{\Delta t}w$ 
yields an approximation of the matrix $B_h$ 
by the matrix $B_h\appr:= Z_h^* [S^{\Delta t}W_h]\appr $, whose
entries are given, for $1\leq i \leq p$ and $1\leq j\leq q$, by:
\begin{eqnarray*}
(B_h\appr )_{ji}&=&
\sup_{x\in X}(z_j(x)+w_i(x)+\Delta t H(x,\frac{\partial w_i}{\partial
    x}))\enspace.
\end{eqnarray*}
Thus, computing $B_h\appr$ requires to solve an optimization problem,
which is nothing but a perturbation of the optimization problem
associated to the computation of $A_h$. 
We may exploit this observation by replacing
$B_h\appr$ by the matrix $B_h\apprd$ with entries
\begin{eqnarray}
(B_h\apprd )_{ji}&=&\<z_j,w_i>+\Delta t \sup_{x\in
  \argmax\{z_j+w_i\}}H(x,\frac{\partial w_i}{\partial x})\enspace,
\label{e-convenient}
\end{eqnarray}
for $1\leq i \leq p$ and $1\leq j\leq q$. Here, 
$\argmax\{z_j+w_i\}$ denotes the set of $x$ such that 
$z_j(x)+w_i(x)=\<z_j,w_i>$. When this set has only one element,
\eqref{e-convenient} yields a convenient approximation
of $B_h$.

Of course, $w_i$ must be differentiable for the approximation~\eqref{stilde}
to make sense. When $w_i$ is non-differentiable, but $z_j$ is differentiable, we may approximate $(B_h)_{ji}$ by
\begin{eqnarray*}
\sup_{x\in X}(z_j(x)+\Delta t H(x,-\frac{\partial z_j}{\partial
    x}) +w_i(x))\enspace,
\end{eqnarray*}
using the dual formula~\eqref{matrixB'}. We may also use the dual 
formula of~\eqref{e-convenient}, where ${\partial w_i}\over {\partial x}$
is replaced by $-{{\partial z_j}\over {\partial x}}$.
\subsection{Comparison with the method of Fleming and McEneaney}
Fleming and McEneaney proposed a max-plus based
method~\cite{mceneaney}, which also uses a space $\calw_h$ generated
by finite elements, $w_1,\ldots,w_p$, together with the linear
formulation~\eqref{exact}.
Their method approaches the value function at time $t$, 
$v^t$, by $W_h \mu^t$, where $W_h=\col(w_i)_{1\leq i\leq p}$ as
above, and $\mu^t$ is defined inductively by
\begin{subequations}\label{algo-fm}
\begin{align}
\mu^0 &= W_h\lres \phi \\
\mu^{t+\Delta t} & =  \big(W_h\backslash (S^{\Delta t}W_h)\big)\mu^{t} \enspace,\label{mut+dt}
\end{align}\end{subequations}
for $t=0,\Delta t,\ldots, T-\Delta t$. This can be compared with the limit
case of our finite element method, in which
the space of test functions $\calz_h$ 
generates the set of all functions. This limit case
corresponds to replacing $Z_h^*$ by the identity operator
in~\eqref{lambdat+dt}, so that 
\begin{align}
\label{e-limit}
\lambda^{t+\Delta t}=W_h\backslash (S^{\Delta t}W_h\lambda^t) \enspace .
\end{align}
\begin{prop}
Let $(\mu^t)$ be the sequence of vectors
defined by the algorithm of Fleming and McEneaney, \eqref{algo-fm};
let $(\lambda^t)$ be the sequence of vectors defined
by the max-plus finite element method, in the limit case~\eqref{e-limit};
and let $v^t$ denote the value function at time $t$. 
Then, 
\[
W_h \mu^t \leq W_h \lambda^t \leq v^t \enspace ,
\quad \mrm{for } t=0,\Delta t, \ldots, T \enspace .
\]
\end{prop}
\begin{proof}[Sketch of proof]
This can be proved by induction, by using the residuation inequality
$W_h\sh S^{\Delta t}W_h\lambda\geq \big(W_h\backslash
(S^{\Delta t}W_h)\big)\lambda$, which holds for all 
vectors $\lambda$, together with the monotonicity
of the operators arising in the construction of $\lambda^t$ and $\mu^t$. 
\end{proof}
An approximation of \eqref{mut+dt} using formulae of the same type as
\eqref{stilde} is also discussed in~\cite{mceneaney99}.
An experimental comparison will appear elsewhere.
\section{Error analysis}
The following general lemma shows that the error of the finite element
method is controlled by the projection errors, 
$\|\projimker{W_h}{}v^t-v^t\|_{\infty}$ and
$\|\projimker{}{Z_h^*}v^t-v^t\|_{\infty}$,  and by the approximation
errors, $\|[S^{\Delta t}w_i]\appr -S^{\Delta t}w_i\|_{\infty}$,
and $|(B_h\apprd)_{ji}-(B_h\appr)_{ji}|$.
\begin{lem}
For $t=0,\Delta t,\cdots,T$,
let $v^t$ be the value function at time $t$, and $v_h^t$
be its approximation given by the effective max-plus finite element method,
implemented with the approximation $B_h\apprd$ of $B_h$,
given by~\eqref{e-convenient}.
We have
\begin{align*}
\|v_h^T- v^T\|_{\infty}&\leq (1+\frac{T}{\Delta
    t})\Big(\sup_{t=0,\Delta t, \ldots,T}(\|\projimker{}{Z_h^*}v^t-v^t\|_{\infty}+\|\projimker{W_h}{}v^t-v^t\|_{\infty})\\
&+\max_{1\leq i\leq p}\|[S^{\Delta t}w_i]\appr -S^{\Delta t}w_i\|_{\infty}+\max_{\begin{subarray}{c}
1\leq j\leq q\\
1\leq i\leq p
\end{subarray}}|(B_h\apprd)_{ji}-(B_h\appr)_{ji}|\Big).
\end{align*}
\end{lem}
The proof of this lemma uses the fact that projectors over max-plus
semimodules are non-expansive in the sup-norm.

To state an error estimate, we make the following assumptions:
\begin{itemize}
\item[-] $(H1)$ The semigroup preserves the 
set of $\frac{1}{c}$-semiconvex functions, for some $c>0$.
\item[-] $(H2)$ $f:X\times U\to \R^n$ is bounded and Lipschitz continuous with
  respect to $x$:
\begin{displaymath}
\begin{array}{cccc}
\exists L_f >0,& \forall x,y\in X, &| f(x,u)-f(y,u)| \leq L_f| x-y |   &\forall u \in U,\\
\exists M_f >0,& \forall x,y\in X,&| f(x,u)|\leq M_f.
\end{array}
\end{displaymath}
\item[-] $(H3)$ $\ell:X\times U\to \R$ is bounded and Lipschitz continuous with respect to $x$:
\begin{displaymath}
\left\{\begin{array}{cc}
| \ell(x,u)-\ell(y,u) |\leq L_l| x-y | & \forall x,y \in X, u \in U,\\
| \ell(x,u) |\leq M_l, & \forall  x,y \in X, u \in U.
\end{array}
\right.
\end{displaymath}
\item[-] $(H4)$ $\phi:X \to \R$ is bounded and Lipschitz continuous:
\begin{displaymath}
\begin{array}{cc}
|\phi(x)-\phi(y)|\leq L_{\phi}| x-y| & \forall x,y\in X.
\end{array}
\end{displaymath} 
\end{itemize}
Recall that a function $f$ is \new{$\frac 1 c$-semiconvex} if
$f(x)+\frac 1{2c}x^2$ is convex. Spaces of semiconvex functions
were already used by Fleming and McEneaney~\cite{mceneaney}.

We shall use the following finite elements.
\begin{defin}[Lipschitz finite elements]
Assume that $X$ is an interval of $\R$.
We call \new{Lipschitz finite element} centered at point $\hat x\in X$,
with constant $A>0$, the function $w(x)=-A|x-\hat x|$.
\end{defin}
\begin{defin}[Quadratic finite elements]
Assume that $X$ is an interval of $\R$.
We call \new{quadratic finite element} centered at point $\hat x\in X$,
with Hessian $\frac 1 c>0$, the function $w(x)=-\frac{1}{2c}(x-\hat x)^2$.
\end{defin}
The family of Lipschitz continuous finite elements of constant $A$ generates,
in the max-plus sense, the semimodule of Lipschitz continuous functions
of Lipschitz constant $A$.
When $X=\R$,
the family of quadratic finite elements with Hessian $\frac 1 c$ generates,
in the max-plus sense, the semimodule of lower-semicontinuous 
$\frac 1 c$-semi-convex functions.
\begin{thm}\label{th-main}
Let $X=[-b,b]\subset \R$. 
We make assumptions \mrm{(H1)}-\mrm{(H4)},
and assume that there exist $L>0$ such that the value function at time
$t$, $v^t$, is $L$-Lipschitz continuous and
$\frac{1}{c}$-semiconvex for all $t>0$,
with the same constant $c$ as in \mrm{(H1)}.
Let us choose quadratic finite elements 
$w_i$ of Hessian $\frac 1 c$, centered at the
points of the regular grid $(\Z \Delta x)\cap [-(b+cL),(b+cL)]$.
Let us choose, as test functions $z_j$,
the Lipschitz finite elements with constant $A\geq L$,
centered at the points of the regular grid
$(\Z \Delta x)\cap [-b,b]$. 
For $t=0,\Delta t,\ldots, T$, let $v_h^t$ be the approximation 
of $v^t$ given by the effective max-plus finite element method,
implemented with the approximation $B_h\apprd$ of $B_h$.
Then, there exists a constant $K>0$ such that,
for $\Delta t$ small enough, 
\[
\|v_h^T-v^T\|_{\infty}\leq K(\sqrt{\Delta t}+\frac{\Delta x}{\Delta t})
\enspace .
\]
\end{thm}
A variant of this theorem, with a stronger assumption, is proved
in~\cite{asma}. We shall give elsewhere the proof of Theorem~\ref{th-main}.
\section{Numerical results}
\begin{exmp}[Linear Quadratic Problem]\label{ex-lq}
We consider the case where $U=\R$, $X=\R$,
\[
\ell(x,u)=-(\frac{a}{2}| x |^2+\frac{| u|^2}{2}),
\quad f(x,u)=u,\mrm{ and } \phi\equiv 0
\enspace .
\]
We obtain $H(x,p)=-\frac{a}{2}| x|^2+\frac{p^2}{2}$. We choose
quadratic finite elements 
$w_i$ and $z_j$ of Hessian $1$, centered at the
points of the regular grid $(\Z \Delta x)\cap [-L,L]$.
We represent in Figure~\ref{quadratique} the solution given by
our algorithm in the case where $T=5$, $\Delta t=\Delta x=0.05$,
$a=0.3$ and $L=10$. The computations were performed using the max-plus
toolbox of Scilab~\cite{toolbox}.
\begin{figure}[htbp]
\begin{center}
\includegraphics[scale=0.4]{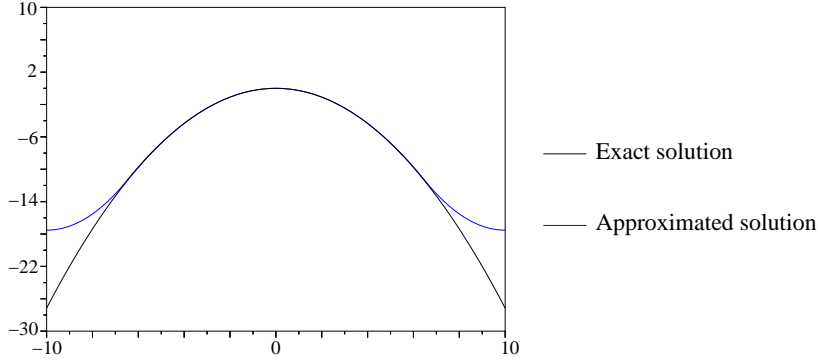}
\caption{Max-plus approximation of a linear quadratic control problem (Example~\ref{ex-lq})}
\label{quadratique}
\end{center}
\end{figure}
\end{exmp}
\begin{exmp}[Distance problem]\label{ex-dist}
We consider the case where $T=1$, $\phi\equiv 0$, $X=[-1,1]$, $U=[-1,1]$,  
\[
\ell(x,u)=\begin{cases}
-1 & \mathrm{if} \quad x \in (-1,1),\\
 0 & \mathrm{if} \quad x \in \{-1,1\} ,
\end{cases}
\quad \mrm{and}\quad
f(x,u)=\begin{cases}
u & \mathrm{if} \quad x \in (-1,1),\\
0 & \mathrm{if} \quad x \in \{-1,1\}.
\end{cases}
\]
Consider first quadratic finite elements 
$w_i$ and $z_j$ of Hessian $\frac{1}{c}$, centered at the
points of the regular grid $(\Z \Delta x)\cap [-1,1]$.
In Figure~\ref{distmauvais}, we represent the solution given by
our algorithm in the case where
$\Delta t=0.05$, $\Delta x=0.0125$ and $c=1.2$. 
Since $\Pi^{Z_h^*}$ is a projector on a subsemimodule
of the $\rminb$-semimodule of $-\frac 1 c$-semiconcave functions,
and since the solution is not $-\frac 1 c$-semiconcave
for any $c$, the error of projection
 $\|\Pi^{Z_h^*}(v^t)-v^t\|_\infty$ does not converge
to zero when $\Delta x$ goes to zero,
which explains the magnitude of the error.
\begin{figure}[htbp]
\begin{center}
\includegraphics[scale=0.4]{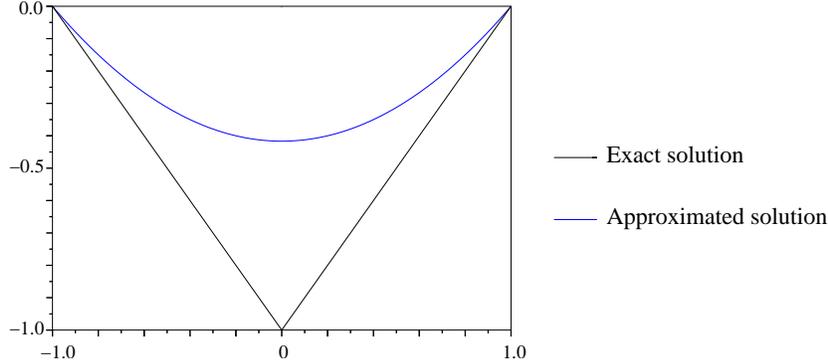}
\caption{A bad choice of test functions for the distance problem (Example~\ref{ex-dist})}
\label{distmauvais}
\end{center}
\end{figure}

To solve this problem, it suffices to replace
the test functions $z_j$ by the Lipschitz finite elements with constant $A\geq 1$, centered at the points of the regular grid
$(\Z \Delta x)\cap [-1,1]$. This is illustrated 
in Figure~\ref{dist} in the case where $\Delta t=0.05$, $\Delta
x=0.0125$, $c=1.2$ and $A=1.1$.
\begin{figure}[htbp]
\begin{center}
\includegraphics[scale=0.4]{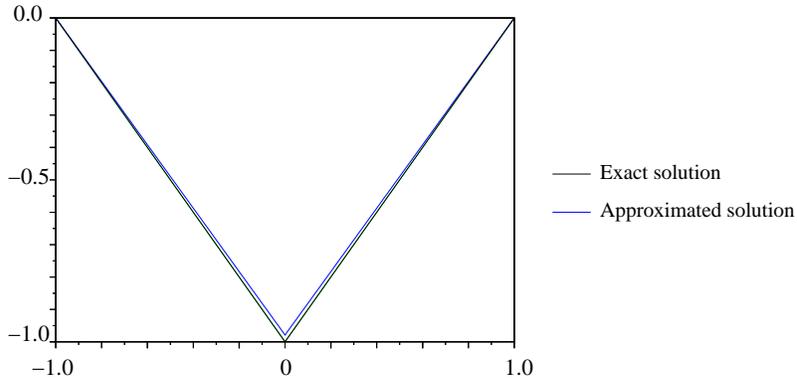}
\caption{A good choice of test functions for the distance problem (Example~\ref{ex-dist})}
\label{dist}
\end{center}
\end{figure}
\end{exmp}
The next two examples are inspired by those proposed by M. Falcone
in \cite{bardi-capuzzo-dolcetta}.
\begin{exmp}\label{ex-falcone1}
We consider the case where $T=1$, $\Phi\equiv 0$, $X=[-1,1]$,
$U=[0,1]$, $\ell(x,u)=x$ and $f(x,u)=-xu$. The optimal choice is to
take $u^*=0$ whenever $x>0$ and to move on the right with maximum
speed ($u^*=1$) whenever $x\leq 0$. For all $t\in [0,T]$, the value function is:
\begin{equation*}
v(x,t)=\begin{cases}
xt & \text{if }  x>0\\
x(1-e^{-t}) & \text{otherwise.}
\end{cases}
\end{equation*}
We choose quadratic finite elements $w_i$ of Hessian $\frac 1 c$ and
Lipschitz finite elements $z_j$ with constant $A\geq 1$. We
represent in Figure~\ref{fal1} the solution given by our
algorithm in the case where $T=1$, $\Delta t=0.05$, $\Delta x=0.02$,
$A=1.3$ and $c=1.4$.
\begin{figure}[htbp]
\begin{center}
\includegraphics[scale=0.4]{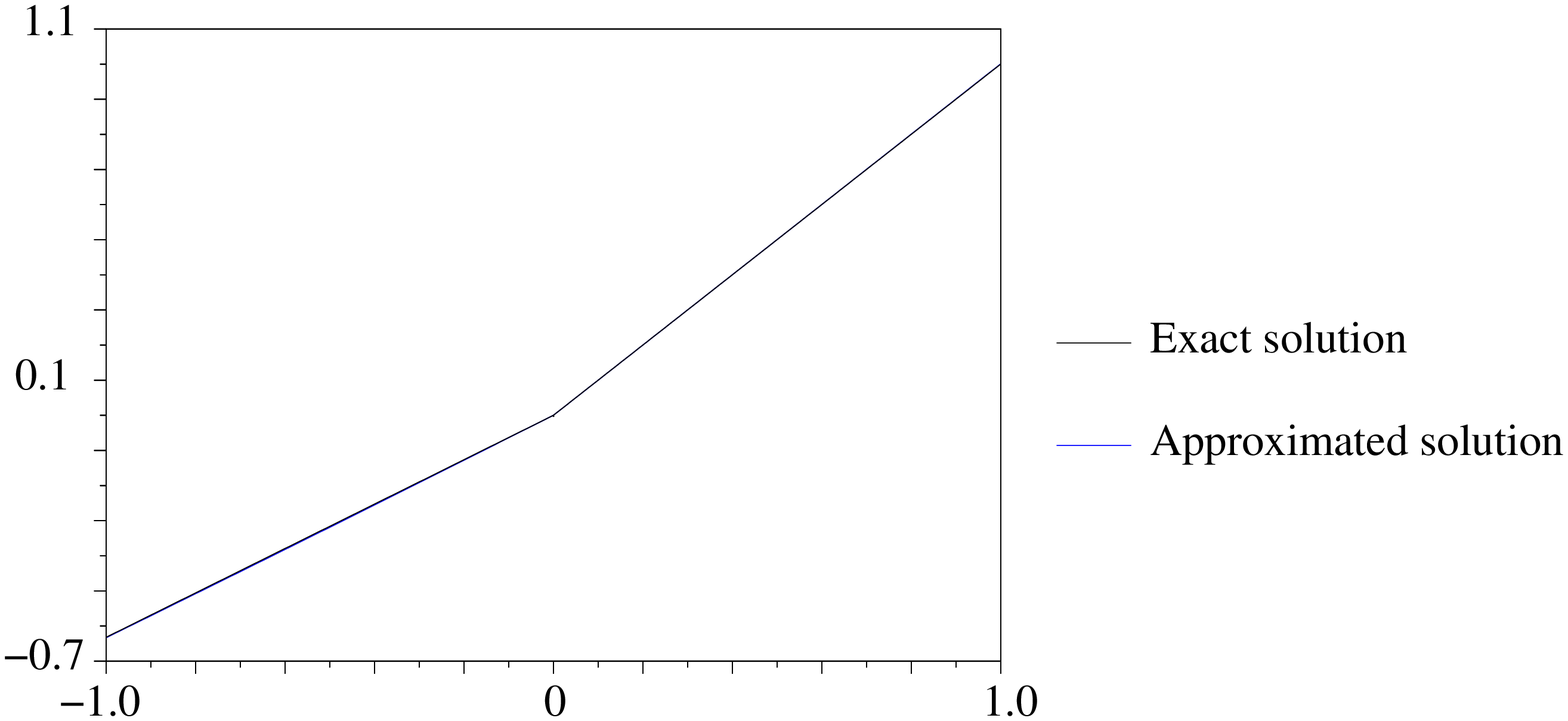}
\caption{Value function and its max-plus approximation (Example~\ref{ex-falcone1})}
\label{fal1}
\end{center}
\end{figure}
\end{exmp}
\begin{exmp}\label{ex-falcone2}
We consider the case where $T=1$, $\Phi\equiv 0$, $X=[-1,1]$,
$U=[-1,1]$, $\ell(x,u)=-3(1-|x|)$ and $f(x,u)=u(1-|x|)$. The optimal
choice is to take $u^*=-1$ whenever $x>0$ and $u^*=1$ whenever
$x<0$. For all $t\in [0,T]$, the value function is:
\begin{equation*}
v(x,t)=-3(1-|x|)(1-e^{-t}).
\end{equation*}
We choose quadratic finite elements $w_i$ of Hessian $\frac 1 c$ and
Lipschitz finite elements $z_j$ with constant $A$. We
represent in Figure~\ref{fal2} the solution given by our
algorithm in the case where $T=1$, $\Delta t=0.05$, $\Delta x=0.02$,
$A=2$ and $c=1.1$.
\begin{figure}[htbp]
\begin{center}
\includegraphics[scale=0.4]{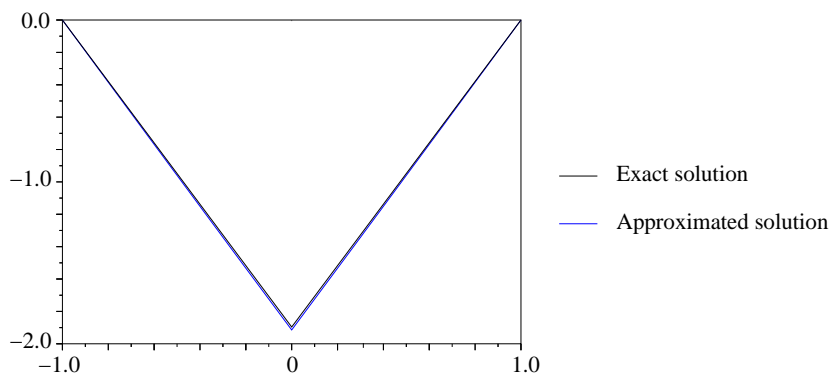}
\caption{Value function and its max-plus approximation (Example~\ref{ex-falcone2})}
\label{fal2}
\end{center}
\end{figure}
\end{exmp}

\bibliography{biblio0304}

\end{document}